\documentclass[10pt]{amsart}

\newlength{\basicwidth} 

\newlength{\shortbasicwidth} 
\newlength{\basicheight} 

\numberwithin{equation}{section}

\begin{document}

\begin{center}
\title{Inequalities for  $1/\bigl(1-\cos(x) \bigr)$ and its derivatives}
\maketitle
\end{center}

\vspace{0.7cm}
\begin{center}
HORST ALZER$^a$ \quad\mbox{and} \quad HENRIK L. PEDERSEN$^b$
\end{center}

\vspace{0.9cm}
\begin{center}
$^a$ Morsbacher Stra\ss {e} 10, 51545 Waldbr\"ol, Germany\\
\emph{Email:} \tt{h.alzer@gmx.de}
\end{center}

\vspace{0.3cm}
\begin{center}
$^b$    Department of Mathematical Sciences, Faculty of Science,\\
 University of Copenhagen, Universitetsparken 5, 2100 Denmark\\
\emph{Email:} \tt{henrikp@math.ku.dk}
\end{center}

\vspace{2.8cm}
{\bf{Abstract.}} We prove that the function $g(x)= 1 / \bigl( 1 - \cos(x)  \bigr)$ is completely monotonic on $(0,\pi]$ and absolutely monotonic on $[\pi, 2\pi)$, and we determine the best possible bounds $\lambda_n$ and $\mu_n$ such that the inequalities
$$
\lambda_n \leq g^{(n)}(x)+g^{(n)}(y)-g^{(n)}(x+y)    \quad (n \geq 0 \,\,\,  \mbox{even})
$$
and
$$
\mu_n \leq g^{(n)}(x+y)-g^{(n)}(x)-g^{(n)}(y)     \quad (n \geq 1 \,\,\,  \mbox{odd})
$$
hold for all $x,y\in (0,\pi)$ with $x+y\leq \pi$.

\vspace{1cm}
{\bf{2020 Mathematics Subject Classification.}} 26A48, 26D05, 11B68.

\vspace{0.08cm}
{\bf{Keywords.}} Completely and absolutely monotonic, trigonometric function, inequality, sub- and superadditive, Bernoulli and Euler numbers, F\'ejer's sine polynomial, $\pi$.

\newpage

\section{Introduction and statement of the results}

{\bf{I.}} \,
A function $f: I\rightarrow \mathbb{R}$, where  $I\subset \mathbb{R}$ is an interval,
is called completely monotonic if $f$ has derivatives of all orders and satisfies
\begin{equation}
(-1)^n f^{(n)}(x)\geq 0 \quad (n=0,1,2,...; \, x\in I).
\end{equation}
If
\begin{equation}
 f^{(n)}(x)\geq 0 \quad (n=0,1,2,...; \, x\in I),
\end{equation}
then $f$ is said to be absolutely monotonic. These functions have remarkable applications in classical analysis, probability theory, numerical analysis and other fields. The main  properties of completely and absolutely monotonic functions can be found in Widder \cite[Chapter IV]{W}.

In the recent past, numerous authors proved that various functions which are defined in terms of well-known special functions, such as gamma and digamma functions, satisfy (1.1) or (1.2). For detailed information on this subject we refer to  Alzer and Berg \cite{AB1, AB2}, Koumandos and Pedersen \cite{KP}, Milovanovi\'c et al. \cite[Chapters 4.2.4, 4.2.5]{MMR}.

In most cases, it was shown that the investigated functions are completely monotonic on  $(0,\infty)$, whereas functions which are completely monotonic only on a bounded interval seem to be less studied.
Here, we present a trigonometric function, defined on $(0, 2\pi)$,  which is completely monotonic on $(0,\pi]$ and absolutely 
monotonic on $[\pi, 2 \pi)$. We apply methods from complex analysis to prove the following result.

\vspace{0.4cm}
\noindent
{\bf{Theorem 1.}} \emph{The function
$$
g(x)=\frac{1}{1-\cos(x)}
$$
is completely monotonic on $(0,\pi]$ and absolutely monotonic on $[\pi, 2\pi)$.}

\vspace{0.4cm}
Our proof yields that the equation
$$
g^{(n)}(x)=0 \quad (0\leq n \in \mathbb{Z}; \, 0<x<2\pi)
$$
holds if and only if $n$ is odd and $x=\pi$. Applications of  Theorem 1 lead to three completely monotonic functions which we present in the following  two corollaries.

\vspace{0.4cm}
{\bf{Corollary 1.}} \emph{The functions
$$
h_1(x)= \exp \bigl( \cot(x) \bigr) \quad (0<x\leq \pi/2)
$$
and
$$
h_2(x)= \exp  \left(  \frac{1}{1+\tan(x)} \right)  \quad (-\pi/4 < x \leq \pi/4)
$$
are completely monotonic.}

\vspace{0.4cm}
A paper published by Lampret \cite{L}, who obtained a closed-form formula for all derivatives of the arctangent function, led us to the function
\begin{equation}
H(x)=\sum_{k=1}^\infty F_k(x) (\cos(x))^k \quad\mbox{with} \quad
\quad
F_k(x)=\sum_{\nu=1}^k \frac{\sin(\nu x)}{\nu}.
\end{equation}
The function $F_k$  is known as F\'ejer's sine polynomial. It plays an important role in the theory of trigonometric polynomials;  see \cite[Chapter 4]{MMR}.

\vspace{0.4cm}
{\bf{Corollary 2.}} \emph{Let $x \in (0,\pi)$. We have
\begin{equation}
H(x)=\frac{\pi/2-x}{1-\cos(x)}.
\end{equation}
In particular, $H$
is completely monotonic on $(0,\pi/2]$.}

\vspace{0.4cm}
{\bf{Remark 1.}} 
Formula (1.4) can be used to find series representations  for $\pi$ and other mathematical constants. For example, from (1.3) and (1.4) with $x=\pi/3$ we obtain 
$$
\pi= \frac{3\sqrt{3}}{2}\sum_{k=1}^\infty \frac{ a_k + b_k}{2^k}
$$
with
$$
a_k= \sum_{\nu=0}^{[(k-1)/3]} \frac{(-1)^{\nu}}{3\nu+1} \quad\mbox{and} \quad
b_k=\sum_{\nu=0}^{[(k-2)/3]} \frac{(-1)^{\nu}}{3\nu+2}.
$$

\vspace{0.4cm}
{\bf{Remark 2.}} Our techniques  to settle Theorem 1 and Corollary 1 can be applied to prove the complete monotonicity of related functions. For example,  we obtain that
$$
x\mapsto \frac{1}{\sin(x)} \quad\mbox{and} \quad x\mapsto \frac{\sin(x)}{1-\cos(x)}
$$
are completely monotonic on $(0,\pi/2)$.

\vspace{0.4cm}
{\bf{II.}} \, Let $I\subset \mathbb{R}$ be an interval. A function $f: I \rightarrow \mathbb{R}$ is called subadditive if
$$
f(x+y)\leq f(x)+f(y) \quad (x,y\in I, x+y\in I).
$$
Functions which satisfy the converse inequality
$$
f(x)+f(x)\leq f(x+y) \quad (x,y\in I, x+y\in I)
$$
are said to be superadditive. Sub- and superadditive functions appear in several fields, like, for example, semigroup theory, number theory and information theory. For more information we refer to Kuczma \cite[Section 16]{Ku}.

Applying Theorem 1 gives for $x,y\in (0,\pi)$ with $x+y\leq \pi$,
\begin{equation}
0< g^{(n)}(x)+g^{(n)}(y)-g^{(n)}(x+y), \quad\mbox{if $n\geq 0$ is even},
\end{equation}
and
\begin{equation}
0< g^{(n)}(x+y)-g^{(n)}(x)-g^{(n)}(y), \quad\mbox{if $n\geq 1$ is odd}.
\end{equation}
This implies that if $n$ is even, then $g^{(n)}$ is (strictly) subadditive on $(0,\pi)$, and if $n\geq 1$ is odd, then $g^{(n)}$ is (strictly) superadditive on $(0,\pi)$. Is it prossible to improve (1.5) and (1.6)? More precisely, we ask whether the lower bound $0$ can be replaced by a positive term which depends only on $n$. The following theorem gives an affirmative answer to this question. It turns out that the best possible lower bounds in (1.5) and (1.6) are given in terms of the classical Bernoulli and Euler numbers, $B_n$ and $E_n$.

\vspace{0.4cm}
{\bf{Theorem 2.}} (i) \emph{Let $n\geq 0$ be an even integer. For all $x,y>0$ with $x+y\leq \pi$, we have
\begin{equation}
\lambda_n \leq g^{(n)}(x)+g^{(n)}(y)-g^{(n)}(x+y)
\end{equation}
with the best possible lower bound}
\begin{equation}
\lambda_n= \frac{2}{n+2} (2^{n+2}-1)^2 |B_{n+2}|.
\end{equation}
(ii)  \emph{Let $n\geq 1$ be an odd integer. For all $x,y>0$ with $x+y\leq \pi$, we have
\begin{equation}
\mu_n \leq g^{(n)}(x+y)-g^{(n)}(x)-g^{(n)}(y)
\end{equation}
with the best possible lower bound}
\begin{equation}
\mu_n= 2 |E_{n+1}|.
\end{equation}

\vspace{0.4cm}
{\bf{Remark 3.}} Since
$$
(-1)^n g^{(n)}(x) \rightarrow \infty \quad\mbox{as} \quad x\rightarrow 0^+,
$$
we conclude that there are no upper bounds for the sums in (1.7) and (1.9) which depend only on $n$.

\vspace{0.4cm}
\section{Proofs}

\emph{Proof of Theorem 1.}
(i) First, we consider the meromorphic function
$$
u(z)=\frac{1}{1-\sin(z)}.
$$
The pole set is $S=\{ z_p=\pi/2+2\pi p \, | \, p\in \mathbb{Z} \}$. Each pole is of order $2$. We have
$$
u(z)= \frac{2}{(z-z_p)^2} +\frac{c_{-1}}{z-z_p} + c_0 + c_1 (z-z_p) + \cdots .
$$
It follows that
$$
c_{-1} = \mbox{Res}(u, z_p) =\lim_{z\to z_p} (z-z_p) \Bigl( u(z)- \frac{2}{(z-z_p)^2}\Bigl)= \lim_{z \to z_p} \frac{(z-z_p)^2 -2(1-\sin(z))}{(z-z_p)(1-\sin(z))}=0.
$$
Let $a\in \mathbb{C} \setminus S$ and
$$
v(z)=\frac{u(z)}{z-a}.
$$
We obtain
$$
\mbox{Res}(v,a)=u(a)=\frac{1}{1-\sin(a)}.
$$
Since, for $z$ close to $z_p$,
\begin{eqnarray}\nonumber
v(z) & = & \Bigl(  \frac{1}{z_p -a}-\frac{1}{(z_p -a)^2} (z-z_p)+ \cdots\Bigr)\Bigl( \frac{2}{(z-z_p)^2}+c_0 + \cdots\Bigr) \\ \nonumber
& = & \frac{2}{z_p-a}\frac{1}{(z-z_p)^2} -\frac{2}{(z_p-a)^2} \frac{1}{z-z_p} + d_0 + \cdots, \\ \nonumber
\end{eqnarray}
we get
$$
\mbox{Res} (v,z_p)=\frac{-2}{(z_p -a)^2}.
$$
Let $m\geq 1$ be an integer and let $R_m = \{z=x+iy \, | \, -m\pi  \leq x,y \leq m\pi\}$. 
 On the 
vertical sides of the boundary of $R_m$ $(z=\pm  m\pi+ i y \pi , \, -m \leq y \leq m)$ we have
$$
1-\sin(z) =1-\sin( \pm m\pi) \cosh( y\pi) -i \cos( \pm m\pi) \sinh( y\pi)=1-i (-1)^m \sinh(y\pi),
$$
$$
|1-\sin(z)|^2 =1+\sinh^2(y\pi).
$$
On the horizontal sides of $R_m$ $(z=x \pi  \pm i m\pi, \, -m \leq x \leq m)$ we get
\begin{eqnarray}\nonumber
|1-\sin(z)|^2 & = & | 1-\sin(x \pi ) \cosh(\pm m \pi) -i \cos(x \pi ) \sinh(\pm m \pi)|^2 \\ \nonumber
& = & 1 +\sin^2(x \pi) \cosh^2(m \pi) -2\sin(x \pi ) \cosh(m \pi) +\cos^2(x \pi )\sinh^2(m  \pi ) \\ \nonumber
& \geq & 1+\sinh^2(m \pi) -2 \sin( x\pi) \cosh(m\pi) \\ \nonumber
& \geq & 1+\sinh^2(m \pi )-2\cosh(m \pi) \\ \nonumber
& = & \bigl(  \cosh(m\pi)-2 \bigr) \cosh(m\pi). \nonumber
\end{eqnarray}
We denote by $\partial_{+} R_m$ the boundary of $R_m$ traversed in the positive direction.
Integration gives
$$
\int_{\partial_{+} R_m}^{} v(z) dz = \pi \int_{-m}^m v( x \pi  -i m \pi) dx
- \pi \int_{-m}^m v(x \pi  + i m \pi) dx
+ i \pi \int_{-m}^m v(m \pi  + i  y \pi ) dy
- i  \pi \int_{-m}^m v(- m \pi  + i y \pi) dy.
$$
Let $m>|a|/\pi$. Then, for $z\in \partial R_m$,
$$
|z-a| \geq |z|-|a| \geq m\pi-|a| >0.
$$
We obtain
\begin{eqnarray}\nonumber
\Big{|} \int_{-m}^m v(\pm m \pi  + i y \pi ) dy\Big{|}
 & \leq & \mbox{sup}_{z\in\partial R_m} \frac{1}{|z-a|} \cdot \int_{-m}^m \frac{dy}{ |1-\sin( \pm m\pi +i y \pi)    |} \\ \nonumber
&  \leq & \frac{1}{m \pi -|a| }  \int_{-m}^m \frac{1}{\cosh(y \pi )}dy \\ \nonumber
&  \leq & \frac{2}{m \pi -|a| }  \int_{0}^\infty \frac{1}{\cosh(y \pi )}dy  \nonumber
\end{eqnarray}
and
$$
\Big{|} \int_{-m}^m v(x \pi    \pm i  m \pi) dx\Big{|}
   \leq  \frac{1}{m\pi -|a|} \frac{2m}{\sqrt{(\cosh(m \pi)-2) \cosh(m \pi)} }
\leq \frac{4}{(m\pi-|a|) \pi}.
$$
It follows that
$$
\lim_{m\to \infty} \int_{\partial_{+} R_m}^{} v(z) dz=0.
$$
Next, we apply the residue theorem. We choose $m$ large enough such that $a$ is in the interior of $R_m$. Then
$$
\frac{1}{2\pi i} \int_{\partial_{+} R_m}^{} v(z) dz 
 =  \frac{1}{1-\sin(a)} + \sum_{-m/2-1/4 < p< m/2 -1/4}^{} \frac{-2}{(z_p -a)^2}. 
$$
We let $m\rightarrow \infty$. Then 
$$
u(a)= \frac{1}{1-\sin(a)} =2 \sum_{p\in \mathbb{Z}}^{} \frac{1}{(z_p -a)^2}.
$$
We differentiate $n$ times and obtain for $a\in \mathbb{C} \setminus{S}$,
$$
u^{(n)}(a)= 2\cdot (n+1)! \sum_{p\in\mathbb{Z}}^{} \frac{1}{(z_p -a)^{n+2}}.
$$
Let $a=t\pi$ with $-1/2 \leq t<1/2$. Then
\begin{equation}
u^{(n)}(a)= \frac{2\cdot (n+1)!}{\pi^{n+2}} \sum_{p\in\mathbb{Z}}^{} \frac{1}{(1/2+2p-t)^{n+2}}.
\end{equation}
If $n$ is even, then $u^{(n)}(a)>0$. Let $n$ be odd. We get
\begin{eqnarray}\nonumber
\sum_{p\in\mathbb{Z}}^{} \frac{1}{(1/2+2p-t)^{n+2}} & = & 
\Bigl( \sum_{p\geq 0} +\sum_{p\leq -1}\Bigr) \frac{1}{(1/2+2p-t)^{n+2}} \\ \nonumber
& = & \sum_{k=0}^\infty \Bigl( \frac{1}{(1/2+2k-t)^{n+2}}-\frac{1}{(1/2+2k+1+t)^{n+2}}\Bigr)\geq 0. \nonumber
\end{eqnarray}
Hence,
$$
u^{(n)}(a) \geq 0 \quad (n=0,1,2,...;\, -\pi/2 \leq  a < \pi/2).
$$

Let $x\in (0,\pi]$. Then $-\pi/2 \leq \pi/2-x< \pi/2$.
Since $g(x)=u(\pi/2-x)$, we obtain 
\begin{equation}
(-1)^n g^{(n)}(x)=u^{(n)}(\pi/2-x)\geq 0  \quad(n=0,1,2,...).
\end{equation}
Thus, $g$ is  completely monotonic  on $(0,\pi]$.

\vspace{0.3cm}
(ii) 
 Let $\pi\leq x<2\pi$. Then $0<2\pi -x \leq \pi$. We have
$$
g(x)=g(2\pi - x).
$$
It follows that
$$
g^{(n)}(x)=  (-1)^n g^{(n)}(2\pi-x)   \geq 0 \quad (n=0,1,2,...).
$$
This means that $g$ is absolutely monotonic on $[\pi, 2\pi)$. $\Box$

\vspace{0.4cm}
{\emph{Proof of Corollary 1.}}  If $y'$ is a completely monotonic function on an interval $I$, then $\exp({-y})$ is also completely monotonic on $I$. This known result can be proved by using induction and the Leibniz rule. 

(i) 
Let $I_1=(0,\pi/2]$ and
$$
y_1(x)=-\cot(x).
$$
Then
$$
y'_1(x)=\frac{2}{1-\cos(2x)}=2 g(2x).
$$
From  Theorem  1 we conclude that $y'_1$ is completely monotonic on $I_1$. 
It follows that $h_1=\exp({-y_1})$ is completely monotonic on $I_1$.

(ii) Let  $I_2 = (-\pi/4, \pi/4]$ and
$$
y_2(x)=\frac{-1}{1+\tan(x)}.
$$
We obtain
$$
y'_2(x)=\frac{1}{1+\sin(2x)} =u(-2x).
$$
Using (2.2) gives for $x\in I_2$,
$$
(-1)^n \bigl( y'_2(x) \bigr)^{(n)} =2^n u^{(n)}(-2x) \geq 0.
$$
This means that $y'_2$ is completely monotonic on $I_2$. Thus also $h_2 = \exp(-y_2)$ is
completely monotonic on $I_2$. $\Box$

\vspace{0.4cm}
\emph{Proof of Corollary 2.}
Lampret \cite{L} proved the formula
$$
\frac{\pi}{2}- |x| =\mbox{sgn}(x) \sum_{k=1}^\infty \frac{\sin(k x)}{k} (\cos(x))^k,
$$
where $x \in \mathbb{R}$ with $0<|x|<\pi$. Let $-1<t<1$. We set $x=\arccos(t)\in (0,\pi)$. Then
$$
\frac{\pi}{2}- \arccos(t) = \sum_{k=1}^\infty \frac{\sin(k\arccos(t))}{k} t^k.
$$
Using the Cauchy product formula gives
$$
\frac{1}{1-t}\cdot \Bigl(\frac{\pi}{2}- \arccos(t) \Bigr)=\sum_{k=0}^\infty t^k \cdot \sum_{k=1}^\infty \frac{\sin(k\arccos(t))}{k} t^k=\sum_{k=1}^\infty F_k(\arccos(t)) t^k.
$$
This leads to
$$
\frac {\pi/2-x}{1-\cos(x)}= \sum_{k=1}^\infty F_k(x) (\cos(x))^k.
$$
The function $x \mapsto \pi/2-x$ is completely monotonic on $(-\infty,\pi/2]$. Since the product of completely monotonic functions is completely monotonic, we conclude that
$$
H(x)  = \Bigl( \frac{\pi}{2}-x\Bigr)\cdot \frac{1}{1-\cos(x)}
$$
is completely monotonic on $(0,\pi/2]$.   $\Box$

\vspace{0.4cm}
{\emph{Proof of Theorem 2.}}  Let $0<x \leq \pi-y<\pi$ and
\begin{equation}
P_n(x,y)= g^{(n)}(x)+g^{(n)}(y)-g^{(n)}(x+y).
\end{equation}

{\underline{Case  1.}} $n \geq 0$ is even.

Since $g^{(n+1)}$ is increasing on $(0,\pi]$, we obtain
$$
\frac{\partial}{\partial x} P_n(x,y)=g^{(n+1)}(x)-g^{(n+1)}(x+y)\leq 0.
$$
Thus,
\begin{equation}
P_n(x,y)\geq P_n(\pi-y,y) = g^{(n)}(\pi-y)+g^{(n)}(y)-g^{(n)}(\pi)=Q_n(y), \quad\mbox{say}.
\end{equation}
The function $Q_n$ is decreasing on $(0,\pi/2]$ and increasing on $[\pi/2, \pi)$, so that we get
\begin{equation}
Q_n(y)\geq Q_n(\pi/2)= 2 g^{(n)}(\pi/2)-g^{(n)}(\pi).
\end{equation}
Using (2.1) and (2.2) gives
\begin{equation}
\frac{\pi^{n+2}}{ 2\cdot (n+1)!} g^{(n)}(\pi/2)
 =  2^{n+2}\sum_{p =0}^{\infty} \frac{1}{(2p+1)^{n+2}} 
\end{equation}
and
\begin{equation}
\frac{\pi^{n+2}}{ 2\cdot (n+1)!} g^{(n)}(\pi) 
 =  2  \sum_{p =0}^{\infty} \frac{1}{(2p+1)^{n+2}}.
\end{equation}

Let $\zeta$ be the Riemann zeta function. We have
\begin{equation}
\sum_{p=0}^{\infty}\frac{1}{(2p+1)^x}= \Bigl(1-\frac{1}{2^x}\Bigr) \zeta(x) \quad (x>1).
\end{equation}
Next, we use (2.6), (2.7), (2.8) and the formula
$$
|B_{2m}| = \frac{2\cdot (2m)!}{(2\pi)^{2m}} \zeta(2m) \quad (m\geq 1).
$$
Then we obtain
\begin{eqnarray}
2 g^{(n)}(\pi/2) - g^{(n)}(\pi) & = &
 \frac{2\cdot (n+1)!}{\pi^{n+2}} \Bigl( 2^{n+3}\sum_{p=0}^\infty \frac{1}{(2p+1)^{n+2}}
- 2  \sum_{p=0}^\infty \frac{1}{(2p+1)^{n+2}} \Bigr)  \\ \nonumber
& = & \frac{(n+1)!}{2^n \pi^{n+2}} (2^{n+2}-1)^2 \zeta(n+2) \\ \nonumber
& = & \frac{2}{n+2}(2^{n+2}-1)^2 |B_{n+2}|. \nonumber
\end{eqnarray}
Finally, we apply (2.3), (2.4), (2.5) and (2.9). Then we obtain (1.7). If we set $x=y=\pi/2$, then equality holds in (1.7). This implies that the lower bound given in (1.8) is sharp.

\vspace{0.2cm}
{\underline{Case 2.}} $n\geq 1$ is odd.

Since $g^{(n+1)}$ is decreasing on $(0,\pi]$, we obtain
$$
\frac{\partial}{\partial x} P_n(x,y)=g^{(n+1)}(x)-g^{(n+1)}(x+y)\geq 0
$$
and
\begin{equation}
P_n(x,y)\leq P_n(\pi-y,y) = g^{(n)}(\pi-y)+g^{(n)}(y)=R_n(y), \quad\mbox{say}.
\end{equation}
The function $R_n$ is increasing on $(0,\pi/2]$ and decreasing on $[\pi/2, \pi)$. It follows that
\begin{equation}
R_n(y)\leq R_n(\pi/2)= 2 g^{(n)}(\pi/2).
\end{equation}
We have
\begin{equation}
-g^{(n)}(\pi/2)= \frac{ 2\cdot (n+1)!}{\pi^{n+2}} \sum_{p \in \mathbb{Z}}^{} \frac{1}{(2p+1/2)^{n+2}}.
\end{equation}
Let $\psi$ be the digamma function. Using
$$
\psi^{(m)}(x) =(-1)^{m+1} m! \sum_{k=0}^\infty \frac{1}{(x+k)^{m+1}} \quad (m\geq 1; x>0)
$$
leads to
\begin{equation}
 \sum_{p\in \mathbb{Z}}^{} \frac{1}{(2p+1/2)^{n+2}}
= \sum_{p=0}^\infty \Bigl(  \frac{1}{(2p+1/2)^{n+2}} -   \frac{1}{(2p+3/2)^{n+2}} \Bigr)
=\frac{-1}{ 2^{n+2} (n+1)!  } \bigl(   \psi^{(n+1)}(1/4) -   \psi^{(n+1)}(3/4) \bigr).
\end{equation}
A result of K\"olbig \cite{K} states that
\begin{equation}
\psi^{(2m)}(1/4) -   \psi^{(2m)}(3/4)= - \pi (2\pi)^{2m}|  E_{2m} |  \quad( m\geq 0).
\end{equation}
Combining (2.12), (2.13) and (2.14) gives
\begin{equation}
-g^{(n)}(\pi/2)= |E_{n+1}|.
\end{equation}
From (2.3), (2.10), (2.11) and (2.15) we conclude that (1.9) holds. If $x=y=\pi/2$, then equality holds in (1.9). Thus,  the lower bound $\mu_n$, as given in (1.10), is sharp. $\Box$.

\end{document}